\documentclass[12pt,twoside]{article}
\usepackage[centertags]{amsmath}
\usepackage{amsfonts}
\usepackage{amssymb}
\usepackage{amsthm}
\usepackage{newlfont}
\newtheorem{Conj}{Conjecture}

\newtheorem{prop}{Proposition}[section]

\begin{document}
	
	{\title{Conjectures in  number theory
			}}

	\author{{\bf  Ahmed ASIMI  }\\ 
		Department of Mathematics, 
		Faculty of Sciences\\
		University Ibnou Zohr,
		B.P. 8106 Agadir, Morroco\\
		\normalsize Laboratory : Laboratoire des Syst\`emes informatique et Vision (LabSiV)\\
		$asimiahmed2013@gmail.com$, $asimiahmed@uiz.ac.ma$ }
	\date{}
	\maketitle
	\begin{abstract}
		Prime numbers, whose properties are important subjects in mathematics, are also fundamental in  computer science notably in IT security, Cryptocurrencies as Bitcoin and Blockchain, cryptography, Code theory notably Error detection codes,  integer factorization, and random number generation.  Finding prime numbers is too an active area of research in mathematics. There are many methods for identifying and generating them and many primality tests which are often complex and expensive in terms of calculation time, and many conjectures and theorems related to prime numbers, such as the prime number theorem, Goldbach's conjecture, and the Riemann hypothesis. My objective in this work is to propose two conjectures :		
	 $1)$ the integer 	
		$1+3.2^{20n}$ is not a prime number for all $n=0$ or $1$ mod $3$, and $2)$ the integer  	
		$1+3.2^{4+20n}$ is not a prime number for all $n=2$  mod $3$. 

	\end{abstract}

\section{ Introduction }
	Prime numbers have fundamental importance in mathematics, and they are used in many fields, including cryptography, number theory, code 
	theory  \cite{LR}, and computer science. There are infinitely many prime numbers, which means there are an infinite amount of them, but their exact distribution in the set of integers has not been fully resolved and is an active subject of mathematical research. They have a central role in the construction, use of finite fields, primitive elemente of finite fields and their arirhmetic, and they are essential for many areas of mathematics and computer science, including cryptography, code theory, and other important applications.
	
	Hence the interest in an in-depth study of their structures, namely the primality of prime numbers, which is our objective in this article in the case where $p=1+3.2^n$. We show that :
	
	\begin{description}
		\item[1)] If $n=3$ mod $4$ then the integer 	
		$1+3.2^{n}$ is not a prime number.
		\item[2)] If $n=2$ mod $3$ then the integer 	
		$1+3.2^{20n}$ is not a prime number for all $n\in \mathbb{N}$.
		\item[3)] We propose two conjectures :
			\begin{description}
			\item[3.1)] The integer 	
			$1+3.2^{20n}$ is not a prime number for all $n=0$ or $1$ mod $3$;
			\item[3.2)] The integer  	
			$1+3.2^{4+20n}$ is not a prime number for all $n=2$  mod $3$.
			
		\end{description}
	\end{description}
		
	\section{ State of the art} 
	
	There are several methods for generating prime numbers, and 	they have many applications in various areas of mathematics and computer science. These applications show how prime numbers are essential to many aspects of modern technology and mathematical research. Their unique property of not being able to be decomposed into products of other numbers makes them valuable for information security and confidentiality.
	
	\textbf{\underline{Methods for generating prime numbers :}}
	
	\begin{description}
		\item[Eratosthenes Sieve: ] This method is one of the oldest for generating prime numbers. It consists of progressively eliminating the multiples of all known prime numbers.
		\item[Fermat's Primality Test: ] Fermat's little theorem \cite{pr1} can be used to test whether a number is prime. However, this test is not perfect and can give false positives for composite numbers, so it is often used as a probabilistic primality test.
		\item[Miller-Rabin Primality Test: ] This is a more sophisticated probabilistic primality test based on Fermat's little theorem. It is used to determine if a number is prime with high probability \cite{pr2}.
		\item[Elliptic curves: ] Elliptic curves are used in cryptography to generate prime numbers which  play a central role in elliptic curves theory by providing the parameters of the curve, defining the order of the points, and ensuring the security of cryptography systems based on these curves. Their  use is complex and essential to ensure the robustness of elliptic curve cryptography systems and other mathematical applications \cite{pr3} and \cite{pr4}. 
		\item[Integer Factoring: ] Another approach to generating prime numbers is to factor a large number into products of two prime numbers. This is used in some encryption algorithms, such as RSA.
		\item[Deterministic primality tests: ] There are deterministic primality tests, such as the AKS primality test, that can determine with certainty whether a number is prime or not. However, these tests are often complex and expensive in terms of calculation time.

\end{description}
	
 \textbf{\underline{Application Areas :}}
	
	\begin{description}
		\item[Cryptography: ] Prime numbers are fundamental in cryptography. They are used in the design of cryptographie hash fuctions, digital signature protocols such as Digital Signature Algorithm (DSA) and Digital Signature Standard (DSS), encryption algorithms such as RSA (Rivest-Shamir-Adleman), where the security of the system relies on the difficulty of factoring a large number into products of prime numbers. Additionally, prime numbers are used to generate cryptographic keys, and also to  solve problems the key exchange protocols, and  the encryption protocols \cite{DH}, \cite{dh} and \cite{NH}.
	\item[Computer security: ] Prime number generation is essential for securing communications and sensitive data over the Internet. They are used to create SSL/TLS certificates for secure websites, to secure protocols such as the Diffie-Hellman protocol \cite{DH}, and also used to generate secure cryptographic keys, protect online communications, and ensure data confidentiality.
	\item[Distributed computing: ] Prime numbers are used in distributed computing projects for scientific research, where volunteers share the computing power of their computers to test the primality of large numbers, thereby contributing to the search for new prime numbers.
	\item[Cryptocurrencies:] Cryptocurrencies, such as blockchain which is a record of all transactions maintained and held by currency holders, rely on mathematical concepts, notably prime numbers, to secure transactions and the creation of new cryptocurrency units.
	\item[Code theory:]  Prime numbers are essential in the theory and practice of error-correcting codes to detect and correct errors in data transmission, particularly in satellite communications and wireless transmissions. Their use in code design helps improve the reliability of digital transmissions by effectively detecting and correcting errors, whether in communications, data storage, or the transmission of digital messages such as Reed-Solomon Codes, Bose-Chaudhuri-Hocquenghem (BCH) and Low-Density Parity-Check (LDPC).
	\item[Number Theory: ] Prime numbers play a central role in number theory. Many unsolved mathematical problems are related to prime numbers, such as the Riemann conjecture, the Goldbach conjecture, the primitive elements of a finite field, Structuring of finite fields, and the prime number theorem.
	\item[Calculation in computer science: ] Prime numbers are used in the optimization of algorithms and data structures, particularly in the design of hash tables and search algorithms.
\item[Random number generation:] Prime numbers are sometimes used to generate high-quality random numbers in computer science and statistical applications. They play also an essential role in the design of Random Number Generators (RNGs), particularly in the context of cryptographically secure random number generation like Cryptographic Security. They are widely used in cryptographic random number generators to ensure data security. The mathematical properties of prime numbers, such as their irreducibility and uniform distribution, contribute to the creation of robust and unpredictable random numbers, which is essential for the security of cryptographic systems \cite{rd1, rd2}.

	\end{description}

	\section{ Our Conjectures}

		The finite fields,  denoted as $\mathbb{F}_{p^{m}}$ where $p$ is a prime number and $m\in \mathbb{N}^{*}$, 	
	 are one of the most beautiful algebraic structures. They are the basis of many algorithmic applications, notably in cryptography, computer science, IT security and code theory whose are currently perfectly dependented. This has led researchers in a natural way to consider methods based on finite fields, elliptic curves and Random number generation in order to construct cryptographic schemes, such as schemes for secure authentication, secret sharing, broadcast encryption,  secure multicast, and cloud computing and IoT by using homomorphic encryption schemes. These structures are perfectly linked to the prime numbers and  the determination of the primitive elements of finite fields. However, there is no feasible and time-resolvable algorithm to test the primality of all integers. 	This innovative scientific work  explores two conjectures in number theory :
	  $1)$ the integer 	
	 $1+3.2^{20n}$ is not a prime number for all $n=0$ or $1$ mod $3$, and $2)$ 	
	 $1+3.2^{4+20n}$ is not a prime number for all $n=2$  mod $3$. 
	 
	 
	 \begin{prop} If $n=3$ mod $4$ then the integer 	
	 	$1+3.2^{n}$ is not a prime number. 
	 \end{prop}
 
 \textbf{Proof.} 
 If $n=3$  mod $4$, then $n=3+4k$, $2^{4}=1$ mod $5$ and  $$\begin{array}[t]{lll}
 1+3.2^{n}	& = & 1+3.2^{3+4k}\\
 	& = & 1+3.2^{3} \mod 5\\
 	& = & 0 \mod 5
 \end{array}$$
  Therefore $5$ devides $1+3.2^{n}$ for all  $n=3$  mod $4$.
  
   \begin{prop} If $n=2$ mod $3$ then the integer 	
  	$1+3.2^{20n}$ is not a prime number for all $n\in \mathbb{N}$. 
  \end{prop}

 \textbf{Proof.} 

If $n=2$ mod $3$, then $n=2+3k$, $2^{3}=1$ mod $7$, $2^{40}=2$ mod $7$ and  $$\begin{array}[t]{lll}
1+3.2^{20n}	& = & 1+3.2^{40}2^{60k}\\
	& = & 1+3.2^{40} \mod 7\\
	& = & 0 \mod 7
\end{array}$$
Therefore $7$ devides $1+3.2^{20n}$ for all  $n=2$  mod $3$.

	\begin{Conj}
	The integer 	
	$1+3.2^{20n}$ is not a prime number for all $n=0$ or $1$ mod $3$. 
		\end{Conj}
	
	 \begin{prop} If $n=0$ or $1$ mod $3$ then the integer 	
		$1+3.2^{4+20n}$ is not a prime number. 
	\end{prop}
	\textbf{Proof.} 
	If $n=0$  mod $3$, then $n=3k$, $2^{3}=1$ mod $7$ and  $$\begin{array}[t]{lll}
		1+3.2^{4+20n}& = & 1+3.2^{4+60k}\\
		& = & 1+3.2^{4} \mod 7\\
		& = & 0 \mod 7
	\end{array}$$
	Therefore $7$ devides $1+3.2^{4+20n}$ for all  $n=0$  mod $3$.
	
	If $n=1$  mod $3$, then $n=1+3k$, $2^{60}=1$ mod $61$, $2^{24}=20$ mod $61$  and  $$\begin{array}[t]{lll}
	 1+3.2^{4+20n}	& = & 1+3.2^{24+60k}\\
		& = & 1+3.2^{24} \mod 61\\
		& = & 1+3.20 \mod 61\\
		& = & 0 \mod 61
	\end{array}$$
	Therefore $61$ devides $1+3.2^{4+20n}$ for all  $n=1$  mod $3$. 
	
	\begin{Conj}
		The integer 	
		$1+3.2^{4+20n}$ is not a prime number for all $n=2$  mod $3$. 
	\end{Conj}

\end{document}